\documentclass[12pt]{amsart}
\usepackage{amssymb}
\usepackage{amsmath}
\usepackage[active]{srcltx}
\usepackage{t1enc}
\usepackage[latin2]{inputenc}
\usepackage{verbatim}
\usepackage{amsmath,amsfonts,amssymb,amsthm}
\usepackage[mathcal]{eucal}
\usepackage{enumerate}
\usepackage[centertags]{amsmath}
\usepackage{graphics}

\newtheorem{theorem}{Theorem}

\newtheorem{corollary}{Corollary}

\newtheorem{proposition}{Proposition}

\begin{document}
\author{Nata Gogolashvili and George Tephnadze}
\title[$T$ means]{On the maximal operators of $T$ means with respect to Walsh-Kaczmarz system}
\address{G. Tephnadze, The University of Georgia, School of science and technology, 77a Merab Kostava St, Tbilisi 0128, Georgia.}
\email{g.tephnadze@ug.edu.ge}
\address{N. Gogolashvili, The University of Georgia, School of science and technology, 77a Merab Kostava St, Tbilisi 0128, Georgia and Department of Computer Science and Computational Engineering, UiT - The Arctic University of Norway, P.O. Box 385, N-8505, Narvik, Norway.}
\email{nata.gogolashvili@gmail.com}
\thanks{The research was supported by Shota Rustaveli National Science	Foundation grant no. FR-19-676.}
\date{}
\maketitle
	
\begin{abstract}
In this paper we prove and discuss some new 
$\left( H_p,L_{p,\infty}\right)$ type inequalities of the maximal operators of $T$ means with monotone coefficients with respect to Walsh-Kaczmarz system. It is also proved that these results are the best possible in a special sense.  As applications, both some well-known and new results are pointed out. In particular, we apply these results to prove a.e. convergence of  such $T$ means. 
\end{abstract}
	
	
	
\section{Introduction}

Concerning definitions and notations used in this introduction we refer to Sections 2. 
	
In 1948, $\breve{\text{S}}$neider \cite{sne} introduced the Walsh-Kaczmarz system and showed that the inequality $\limsup_{n\rightarrow \infty}D_{n}^{\kappa }(x)/\log n\geq C>0$ holds a.e. In 1974 Schipp \cite{Sch2} and Young \cite{Y} proved that the Walsh-Kaczmarz system is a convergence system. 

In \cite{G-N-Complex} Goginava and Nagy proved that the maximal operator $R^{\ast,\kappa }$ of Riesz means with respect to Walsh-Kaczmarz system is bounded from the Hardy space $H_{1/2}$ to the space $L_{1/2,\infty}$, but is not bounded from the Hardy space $H_p$ to the space $L_p,$ when $0<p\leq 1/2.$ In \cite{tep0} it was proved that there exists a martingale $f\in H_{p},$ $(0<p\leq 1),$ such that the maximal operator  $L^{\ast,\kappa }$ of N\"orlund logarithmic means  with respect to Walsh-Kaczmarz system is not bounded in the Lebesgue space $L_{p}.$ 
The Logarithmic means with respect to the Walsh and Vilenkin systems systems were studied by Blahota and Gát \cite{bg}, Gát \cite{Ga1},  Lukkassen, Persson, Tephnadze and Tutberidze \cite{lptw} (see also \cite{GNT1}, \cite{ptw2}, \cite{pttw}, \cite{tep4}, \cite{TT1} and \cite{tut2}), Simon \cite{Si1}.

In 1981, Skvortsov \cite{Sk1} showed that the Fejér means with respect to the Walsh-Kaczmarz system converge uniformly to $f$ for any continuous functions $f$. Gát \cite{gat} proved that, for any integrable functions, the Fejér means with respect to the Walsh-Kaczmarz system converges almost everywhere to the function. He showed that the maximal operator $\sigma^{\ast ,\kappa }$ of \ Walsh-Kaczmarz-Fejér means is of weak type $(1,1)$ and of type $(p,p)$ for all $1<p\leq \infty $. Gát's result was generalized by Simon \cite{S2}, who showed that the maximal operator $\sigma ^{\ast ,\kappa }$ is of type $(H_{p},L_{p})$ for $p>1/2$. In the endpoint case $p=1/2$ Goginava \cite{Gog-PM} (see also  \cite{PTT,tep1,tep2,tep3}) proved that maximal operator 
$\sigma^{\ast,\kappa }$ of Walsh-Kaczmarz-Fejér means is not of type $(H_{1/2},L_{1/2})$. Weisz \cite{We5} showed that the following is true:
	
\textbf{Theorem W1.} The maximal operator $\sigma ^{\ast,\kappa}$ of  Walsh-Kaczmarz-Fejér means is bounded from the Hardy space $H_{1/2}$ to the space $L_{1/2,\infty }$.
	
Gát and Goginava \cite{gago} proved that the  maximal operator 
$\sigma^{\alpha,\ast,\kappa}$ of $\left(C,\alpha \right) $ $\left( 0<\alpha <1\right) $ means with respect Walsh-Kaczmarz system is bounded from the Hardy space $H_{1/\left( 1+\alpha \right) }$ to the space $L_{1/\left( 1+\alpha \right),\infty }$. Goginava and Nagy \cite{gona} proved that $\sigma ^{\alpha ,\ast ,\kappa }$ is not bounded from the Hardy space $H_{1/\left( 1+\alpha \right) }$ to the space $L_{1/\left( 1+\alpha \right) }$. Móricz and Siddiqi \cite{mor} investigated the approximation properties of some special Nörlund means  of $L_p$ function in norm. These means in the martingale Hardy spaces were discussed in  Blahota, Tephnadze \cite{bt1,bt2}. In \cite{ptw} and \cite{tep5}  (see also \cite{BPTW,bpt1,mpt}) it was proved some $(H_p,L_p)$-type inequalities for the maximal operators of Nörlund means, with respect to Walsh-Kaczmarz and Vilenkin systems, when $0<p\leq 1.$ In the two dimensional case approximation properties of Nörlund and Cesáro means were considered by Nagy \cite{na,n,nagy} and by Nagy and Tephnadze  \cite{nt1,nt2,nt3,nt4,nt5}. Some boundedness results of so called $T$, $\Theta$ and $\theta$ means in the Lebesgue and martingale Hardy spaces can be found in Blahota and Nagy \cite{BN}, in  Blahota, Nagy and Tephnadze \cite{BNT}, Tutberidze \cite{tut1} and Weisz \cite{WeTh2,WeTh3,WeTh4,WeTh5}.
	
The main aim of this paper is to investigate $(H_{p},L_{p,\infty })$-type inequalities for the maximal operators of $T$ means with monotone coefficients of the one-dimensional Kaczmarz-Fourier series.
	
This paper is organized as follows: in order not to disturb our discussions	later on some definitions and notations are presented in Section 2. The main results and some of its consequences can be found in Section 3. For the proofs of the main results we need some auxiliary results of independent interest. Also these results are presented in Section 3. The detailed proofs are given in Section 4.
	
	\section{Definitions and Notations}
	
Now, we give a brief introduction to the theory of dyadic analysis \cite{SWSP}. Let $\mathbf{N}_{+}$ denote the set of positive integers, $\mathbf{N:=N}_{+}\cup \{0\}.$
	
Denote ${\mathbb{Z}}_{2}$ the discrete cyclic group of order 2, that is ${\mathbb{Z}}_{2}=\{0,1\},$ where the group operation is the modulo 2 addition and every subset is open. The Haar measure on ${\mathbb{Z}}_{2}$ is given such that the measure of a singleton is 1/2. Let $G$ be the complete direct product of the countable infinite copies of the compact groups ${\mathbb{Z}}_{2}.$ The elements of $G$ are of the form 
\begin{equation*}
x=\left( x_{0},x_{1},...,x_{k},...\right) ,\text{ \ \ \ }x_{k}=0\vee 1,\text{\ }\left( k\in \mathbf{N}\right) .
\end{equation*}
	
The group operation on $G$ is the coordinate-wise addition, the measure (denoted by $\mu $) and the topology are the product measure and topology.The compact Abelian group $G$ is called the Walsh group. A base for the neighborhoods of $G$ can be given in the following way: 
\begin{equation*}
I_{0}\left( x\right) :=G,\quad I_{n}\left( x\right) :=I_{n}\left(
x_{0},...,x_{n-1}\right) :=\left\{ y\in G:\,y=\left(
x_{0},...,x_{n-1},y_{n},y_{n+1},...\right) \right\} ,
\end{equation*}%
$\left( x\in G,n\in \mathbf{N}\right) .$ These sets are called dyadic  intervals. Denote by $0=\left( 0:i\in \mathbf{N}\right) \in G$ the null element of $G.$ Let $I_{n}:=I_{n}\left( 0\right) ,$  $\overline{I_{n}}:=G\backslash I_{n}\,\left( n\in \mathbf{N}\right).$ Set $e_{n}:=\left(0,...,0,1,0,...\right) \in G,$ the $n$-th coordinate of which is 1 and the rest are zeros $\left( n\in \mathbf{N}\right).$
	
For $k\in \mathbf{N}$ and $x\in G$ let us denote the $k$-th Rademacher function, by 
\begin{equation*}
r_{k}\left( x\right) :=\left( -1\right) ^{x_{k}}.
\end{equation*}
	
Now, define the Walsh system $w:=(w_{n}:n\in \mathbf{N})$ on $G$ as: 
\begin{equation*}
w_{n}(x):=\overset{\infty}{\underset{k=0}{\Pi}}r_k^{n_k}\left(x\right) =r_{\left\vert n\right\vert}\left(x\right) \left(-1\right) ^{\underset{k=0}{\overset{\left\vert n\right\vert-1}{\sum}} n_{k}x_{k}}\text{\qquad}\left(n\in\mathbf{N}\right).
\end{equation*}
	
If $n\in \mathbf{N}$, then $n=\sum\limits_{i=0}^{\infty }n_{i}2^{i}$ can be written, where $n_{i}\in \{0,1\}\quad \left( i\in \mathbf{N}\right) $, i. e. $n$ is expressed in the number system of base 2. Denote $\left\vert n\right\vert :=\max \{j\in \mathbf{N}: \ n_{j}\neq 0\},$ that is $2^{\left\vert n\right\vert }\leq n<2^{\left\vert n\right\vert +1}.$
	
The Walsh-Kaczmarz functions are defined by
\begin{equation*}
\kappa _{n}\left( x\right) :=r_{\left\vert n\right\vert }\left( x\right)
\prod\limits_{k=0}^{\left\vert n\right\vert -1}\left( r_{\left\vert
n\right\vert -1-k}\left( x\right) \right) ^{n_{k}}=r_{\left\vert
n\right\vert }\left( x\right) \left( -1\right)^{\sum\limits_{k=0} ^{\left\vert n\right\vert -1}n_{k}x_{_{\left\vert n\right\vert -1-k}}}.
\end{equation*}
	
The Dirichlet kernels are defined by
\begin{equation*}
D_{0}:=0\text{, \ \ \ \ }D_{n}^{\psi }:=\sum_{i=0\text{ }}^{n-1}\psi _{i},\text{ }\left( \psi =w\text{ or }\psi =\kappa \right).
\end{equation*}
	
The $2^{n}$-th Dirichlet kernels have a closed form (see e.g. \cite{SWSP})
\begin{equation} \label{Dir}
D_{2^{n}}^{w}\left( x\right) =D_{2^{n}}\left( x\right) =D_{2^{n}}^{\kappa}\left( x\right) \text{\thinspace }=\left\{ 
\begin{array}{ll}
2^{n}\text{ \ \ \ \ \ } & x\in I_{n}, \\ 
0 & x\notin I_{n}.
\end{array}
\right.  
\end{equation}
	
The norm (or quasi-norm) of the spaces $L_{p}(G)$ and $L_{p,\infty }\left(G\right)$ are respectively defined by
	
\begin{equation*}
\left\Vert f\right\Vert _{p}^{p}:=\int_{G}\left\vert f\right\vert ^{p}d\mu ,\text{ \ \ }\left\Vert f\right\Vert _{L_{p,\infty }(G)}^{p}:=\underset{\lambda >0}{\sup }\lambda ^{p}\mu \left( f>\lambda \right) ,\text{ \ \ \ }\left( 0<p<\infty \right).
\end{equation*}
	
The partial sums with respect to Walsh and Walsh-Kaczmarz series are defined as follows: 
\begin{equation*}
S_{M}^{\psi }f:=\sum\limits_{i=0}^{M-1}\widehat{f}\left( i\right) \psi _{i},\text{ \ \ }\left( \psi =w\text{ or }\psi =\kappa \right).
\end{equation*}
	
Let $\{q_{k}:k\geq 0\}$ be a sequence of nonnegative numbers. The $n$-th Nörlund and  $T$ means for a Fourier series of $f$  are respectively defined by
\begin{equation*}
t^{\psi }_{n}f=\frac{1}{Q_{n}}\overset{n}{\underset{k=1}{\sum }} q_{n-k}S^{\psi }_{k}f, \text{ \ \ }\left( \psi =w\text{ or }\psi =\kappa \right) 
\end{equation*}
and
\begin{equation} \label{nor}
T^{\psi }_nf:=\frac{1}{Q_n}\overset{n-1}{\underset{k=0}{\sum }} q_kS^{\psi}_kf, \text{ \ \ }\left(\psi=w\text{ or }\psi =\kappa\right).
\end{equation}
where $Q_{n}:=\sum_{k=0}^{n-1}q_{k}.$ It is obvious that 
\begin{equation*}
T^{\psi }_nf\left(x\right)=\underset{G_m}{\int}f\left(t\right)F^{\psi }_n\left(x-t\right) d\mu\left(t\right)
\end{equation*}%
where $	F^{\psi }_n:=\frac{1}{Q_n}\overset{n-1}{\underset{k=0}{\sum }}q_{k}D^{\psi }_k$
is called $T$ kernel.
	
We always assume that $\{q_k:k\geq 0\}$ be a sequence of nonnegative numbers and $q_0>0.$ Then the summability method (\ref{nor}) generated by $\{q_k:k\geq 0\}$ is regular if and only if $	\lim_{n\rightarrow\infty}Q_n=\infty.$
	
Let consider some class of $T$ means with monotone and bounded sequence
$\{q_k:k\in \mathbb{N}\}$, such that
\begin{equation*}
q:=\lim_{n\rightarrow\infty}q_n>c>0.
\end{equation*}
	
Then, it easy to check that
\begin{equation} \label{monotone1}
\frac{q_{n-1}}{Q_n}=O\left(\frac{1}{n}\right),\text{ \  as \ }n\rightarrow \infty.
\end{equation}
	
The $n$-th Fejér means of a function $f$ \ is given by
\begin{equation*}
\sigma _{n}^{\psi }f:=\frac{1}{n}\sum_{k=0}^{n-1}S_{k}^{\psi }f,\text{ \ \ \ }\left( \psi =w\text{ or }\psi =\kappa \right) .
\end{equation*}
	
Fejér kernel is defined in the usual manner
\begin{equation*}
K_{n}^{\psi }:=\frac{1}{n}\overset{n}{\underset{k=1}{\sum }}D_{k}^{\psi },\text{ \ \ \ }\left( \psi =w\text{ or }\psi =\kappa \right) .
\end{equation*}

If we invoke Abel transformation to obtain the following identities, which are very important for the investigations of $T$ summability:
\begin{eqnarray} \label{2b}
Q_n&:=&\overset{n-1}{\underset{j=0}{\sum}}q_j =\overset{n-2}{\underset{j=0}{\sum}}\left(q_{j}-q_{j+1}\right) j+q_{n-1}{(n-1)}
\end{eqnarray}
and
\begin{equation}\label{2c}
F_n^{\psi }=\frac{1}{Q_n}\left(\overset{n-2}{\underset{j=0}{\sum}}\left(q_j-q_{j+1}\right) jK^{\psi }_j+q_{n-1}(n-1)K^{\psi }_{n-1}\right).
\end{equation}
	
The $\left( C,\alpha \right) $-means are defined as
\begin{equation*}
\sigma _{n}^{\alpha ,\psi }f=\frac{1}{A_{n}^{\alpha }}\overset{n}{\underset{k=1}{\sum }}A_{n-k}^{\alpha -1}S_{k}^{\psi }f,\text{ \ \ \ \ }\left( \psi =w
\text{ or }\psi =\kappa \right),
\end{equation*}%
where
\begin{equation*}\label{1d}
A_{0}^{\alpha }=0,\text{ \qquad }A_{n}^{\alpha }=\frac{\left( \alpha+1\right) ...\left( \alpha +n\right) }{n!},~~\alpha \neq -1,-2,...
\end{equation*}
	
It is known that%
\begin{equation*} \label{2d}
A_{n}^{\alpha }\sim n^{\alpha },\text{ \ }A_{n}^{\alpha }-A_{n-1}^{\alpha
}=A_{n}^{\alpha -1},\text{ }\overset{n}{\underset{k=1}{\sum }}%
A_{n-k}^{\alpha -1}=A_{n}^{\alpha }.  
\end{equation*}
	
We also consider "inverse" $\left(C,\alpha\right)$-means, which is example of $T$-means:
\begin{equation*}
U_n^{\alpha,\psi}f:=\frac{1}{A_n^{\alpha}}\overset{n-1}{\underset{k=0}{\sum}}A_{k}^{\alpha-1}S_k^\psi f, ,\text{ \ \ }\left( \psi =w
\text{ or }\psi =\kappa \right).
\end{equation*}
	
Let $\beta_n^{\alpha}$ denote the $T$ mean, where
\begin{equation*}
\left\{q_0=1,\qquad q_k=k^{\alpha-1}:k\in \mathbb{N}_+\right\} ,
\end{equation*}
that is 
\begin{equation*}
V_n^{\alpha, ,\psi}f:=\frac{1}{Q_n}\overset{n}{\underset{k=1}{\sum }}k^{\alpha-1}S_k^\psi f,\text{ \ \ }\left( \psi =w\text{ or }\psi =\kappa \right) \qquad 0<\alpha<1.
\end{equation*}
	
The $n$-th Riesz`s logarithmic mean $R^{\psi }_{n}$ and Nörlund logarithmic mean $L^{\psi }_{n}$ are defined by
	
\begin{equation*}
R_{n}^{\psi }f:=\frac{1}{l_{n}}\overset{n-1}{\underset{k=0}{\sum }}\frac{S_{k}^{\psi }f}{k},\text{ \ }L_{n}^{\psi }f:=\frac{1}{l_{n}}\overset{n-1}{\underset{k=1}{\sum }}\frac{S_{k}^{\psi }f}{n-k},\text{ \ \ }\left( \psi =w\text{ or }\psi =\kappa \right) .
\end{equation*}
respectively, where
$l_{n}:=\sum_{k=1}^{n-1}1/k.$
	
Up to now we have considered $T$ mean in the case when the sequence $\{q_k:k\in\mathbb{N}\}$ is bounded but now we consider $T$ summabilities with unbounded sequence $\{q_k:k\in\mathbb{N}\}$. Let $\alpha\in\mathbb{R}_+,\ \ \beta\in\mathbb{N}_+$ and
\begin{equation*}
\log^{(\beta)}x:=\overset{\beta-\text{times}}{\overbrace{\log ...\log}}x.
\end{equation*}
	
If we define the sequence $\{q_k:k\in \mathbb{N}\}$ by
\begin{equation*}
\left\{q_0=0\text{ \  and \ }q_k=\log^{\left(\beta \right)}k^{\alpha
}:k\in\mathbb{N}_+\right\},
\end{equation*}%
then we get the class of $T$ means with non-decreasing coefficients:
\begin{equation*}
B_n^{\alpha,\beta,\psi }f:=\frac{1}{Q_n}
\sum_{k=1}^{n}\log^{\left(\beta\right)}k^{\alpha}S_k^\psi f,\text{ \ \ }\left( \psi =w\text{ or }\psi =\kappa \right).
\end{equation*}%
	
We note that $B_n^{\alpha,\beta}$ are well-defined for every $n \in \mathbb{N}$
\begin{equation*}
B_n^{\alpha,\beta,\psi}f=\sum_{k=0}^{n-1}\frac{\log^{\left(\beta\right)}k^{\alpha }}{Q_n}S_k^\psi f, \ \ \ \ \text{ \ \ }\left( \psi =w \text{ or }\psi =\kappa \right).
\end{equation*}
	
It is obvious that $\frac{n}{2}\log^{\left(\beta \right)}\frac{n^{\alpha }}{2^{\alpha }}\leq Q_n\leq n\log^{\left(\beta\right)}n^{\alpha}.$ It follows that
\begin{eqnarray} \label{node00}
\frac{q_{n-1}}{Q_n}\leq\frac{c\log^{\left(\beta\right)}\left(n-1\right)^{\alpha}}{n\log^{\left(\beta\right) } n^{\alpha}}=O\left(\frac{1}{n}\right) \rightarrow 0,\text{ \ as \ } n\rightarrow \infty.
\end{eqnarray}
	
Let us define maximal operator of $T$ means by 
\begin{eqnarray*}
T^{\ast,\psi}f:=\sup_{n\in\mathbb{N}}\left\vert T_n^\psi f\right\vert, \ \ \ \text{ \ \ }\left( \psi =w \text{ or }\psi =\kappa \right).
\end{eqnarray*}
	
The well-known example of maximal operator of $T$ means are Fejer and Riesz logarithmic means which are defined by: 
\begin{eqnarray*}
\sigma ^{\ast ,\psi}f:=\sup_{n\in \mathbb{N}}\left\vert \sigma _{n}^{\psi }f\right\vert, \ \ \ \ \ \ 
R^{\ast,\psi}f:=\sup_{n\in\mathbb{N}}\left\vert R_{n}^\psi f\right\vert \ \ \ \ \ \ \text{ \ \ }\left( \psi =w \text{ or }\psi =\kappa \right)
\end{eqnarray*}	
respectively. We also define some new maximal operators:
\begin{eqnarray*}
U^{\alpha,\beta,\ast,\psi}f:=\sup_{n\in\mathbb{N}}\left\vert U_n^{\alpha ,\beta, \psi }f\right\vert, \  
V^{\alpha,\ast,\psi }f:=\sup_{n\in\mathbb{N}}\left\vert V_{n}^{\alpha,\psi}f\right\vert, \ 
B^{\alpha,\ast,\psi }f:=\sup_{n\in\mathbb{N}}\left\vert B_{n}^{\alpha,\psi}f\right\vert \  \left( \psi =w \text{ or }\psi =\kappa \right).
\end{eqnarray*}
For the martingale $f$ we consider the following maximal operators of some Norlund means:
\begin{equation*}
t^{\ast ,\psi }f:=\sup_{n\in \mathbb{N}}\left\vert t_{n}^{\psi }f\right\vert  , \text{ \  }\sigma ^{\alpha,\ast,\psi }f:=\sup_{n\in \mathbb{N}}\left\vert\sigma_{n}^{\alpha,\psi}f\right\vert,  \text{ \  }
L^{\ast ,\psi}f:=\sup_{n\in \mathbb{N}	}\left\vert L_{n}^{\psi }f\right\vert,\text{ \  }(\psi =w\text{ or \ }\psi =\kappa).
\end{equation*}
	
The $\sigma $-algebra generated by the dyadic intervals of measure $2^{-k}$ will be denoted by $F_{k}$ $\left( k\in \mathbf{N}\right) .$ Denote by $f=\left( f^{\left( n\right) },n\in \mathbf{N}\right) $ a martingale with respect to $\left( F_{n},n\in \mathbf{N}\right) $ (for details see, e. g. \cite{we2,we3}). The maximal function of a martingale $f$ is defined by 
$f^{\ast }=\sup\limits_{n\in \mathbf{N}}\left\vert f^{\left( n\right) }\right\vert .$
	
In case $f\in L_{1}\left( G\right) $, the maximal function can also be given by 
\begin{equation*}
f^{\ast }\left( x\right) =\sup\limits_{n\in \mathbf{N}}\frac{1}{\mu \left(I_{n}(x)\right) }\left\vert \int\limits_{I_{n}(x)}f\left( u\right) d\mu\left(u\right)\right\vert ,\ \ x\in G.
\end{equation*}
For $0<p<\infty $ the Hardy martingale space $H_{p}(G)$ consists of all martingales for which
\begin{equation*}
\left\| f\right\| _{H_{p}}:=\left\| f^{*}\right\| _{p}<\infty .
\end{equation*}
	
If $f\in L_{1}\left( G\right) ,$ then it is easy to show that the sequence $\left( S^{\psi}_{2^{n}}f:n\in \mathbf{N}\right) $ is a martingale.
	
The Walsh-Fourier and Walsh-Kaczmarz-Fourier coefficients of $f\in L_{1}\left( G\right) $ are the same as the ones of the martingale $\left( S^{\psi}_{2^{n}}f:n\in \mathbf{N}\right) $ obtained from $f$.
	
If $f$  is a martingale, then the Walsh-Fourier and Walsh-Kaczmarz-Fourier coefficients must be defined in a little bit different way: 
\begin{equation*}
\widehat{f}^{\psi }\left( i\right) =\lim\limits_{n\rightarrow \infty}\int\limits_{G}f^{\left( n\right) }\psi _{i}d\mu ,\text{ \ \ }\left( \psi=w,\text{ or }\psi =\kappa \right) .
\end{equation*}
	
A bounded measurable function $a$ is p-atom, if there exists an interval $I$, such that
\begin{equation*}
\int_{I}ad\mu =0,\text{ \ \ }\left\Vert a\right\Vert _{\infty }\leq \mu
\left( I\right) ^{-1/p},\text{ \ \ supp}\left( a\right) \subset I.
\end{equation*}
	
Weisz proved that Hardy spaces $H_p$ have atomic characterization. In particular the following is true:
\begin{proposition} \cite{we2} \label{prop1} A martingale $f=\left( f^{\left( n\right) },n\in\mathbb{N}	\right) $ is in $H_{p}\left( 0<p\leq 1\right) $ if and only if there exists sequence $\left( a_{k},k\in	\mathbb{N}\right) $ of p-atoms and a sequence $\left( \mu _{k},k\in	\mathbb{N} \right) ,$ of real numbers, such that, for every $n\in\mathbb{N},$
\begin{equation} \label{1}
\qquad \sum_{k=0}^{\infty }\mu _{k}S^{\psi}_{2^n}a_{k}=f^{\left( n\right) },\ \ \ \ \ \ \ \sum_{k=0}^{\infty}\left\vert\mu _k\right\vert ^{p}<\infty, \ \ \ (\psi =w\text{ or \ }\psi =\kappa).
\end{equation}
Moreover,
\begin{equation*}
\left\Vert f\right\Vert _{H_{p}}\backsim \inf \left( \sum_{k=0}^{\infty
}\left\vert \mu _{k}\right\vert ^{p}\right) ^{1/p},
\end{equation*}
where the infimum is taken over all decomposition of $f$ of	the form (\ref{1}).
\end{proposition}
	
\begin{proposition} \label{prop2}\cite{we2,we3} 
Suppose that an operator $T$ \ is $\sigma $-linear and for some $0<p<1$ 
\begin{equation*}
\left\Vert Tf\right\Vert _{L_{p,\infty}}\leq c_{p}\left\Vert f\right\Vert_{H_{p}},
\end{equation*}
then $T$ is of weak	type-(1,1):
\begin{equation*}
\left\Vert Tf\right\Vert _{L_{1,\infty}}\leq c\left\Vert f\right\Vert_{1}.
\end{equation*}
\end{proposition}
	
\section{Results}
	
\begin{center}
\textbf{Main results and some of its consequences}
\end{center}
We state our main result concerning the maximal operator of the summation method (\ref{nor}), which we also show is in a sense sharp.
	
\begin{theorem}\label{theorem1}
a) The maximal operator $T^{\ast,\kappa }$ of \ the summability method (\ref{nor}) with non-increasing sequence $\{q_{k}:k\geq 0\},$ is bounded from the Hardy space $H_{1/2}$ to the space $L_{1/2,\infty}.$
		
The statement in a) is sharp in the following sense:
\bigskip b) Let $0<p<1/2$ and\ $\{q_{k}:k\geq 0\}$ is non-decreasing sequence, satisfying the condition
\begin{equation} \label{cond1}
\frac{q_{n+1}}{Q_{n+2}}\geq \frac{c}{n},\text{ \ \ }\left( c\geq 1\right) .
\end{equation}
Then there exists a martingale $f\in H_{p},$ such that
\begin{equation*}
\underset{n\in	\mathbb{N}}{\sup }\left\Vert T^\kappa_{n}f\right\Vert _{L_{p,\infty}}=\infty.
\end{equation*}
\end{theorem}
	
A number of special cases of our results are of particular interest and give both well-known and new information. We just give the following examples of such $T$ means with non-increasing sequence $\{q_{k}:k\geq 0\}:$
\begin{corollary}\label{cor0}
The maximal operators $U^{\alpha,\beta,*,\kappa}$, $V^{\alpha, *,\kappa}$ and $R^{*,\kappa}$ are  bounded from the Hardy space $H_{1/2}$ to the space $L_{1/2,\infty}$ but are not bounded from $H_{p}$ to the space $L_{p,\infty},$ when $0<p<1/2.$
\end{corollary}
	
\begin{corollary}\label{cor2}
Let $f\in L_{1}$ and $T^\kappa_{n}$ be the $T$ means with non-increasing sequence $\{q_{k}:k\geq 0\}$. Then $T^\kappa_{n}f\rightarrow f,\text{ \ \ a.e., \ \ as \ }n\rightarrow \infty.$
\end{corollary}
	
\begin{corollary}\label{cor1}Let $f\in L_1$. Then
\begin{eqnarray*}
R^\kappa_{n}f &\rightarrow &f,\text{ \ \ \ a.e., \ \ \ \ as \ }n\rightarrow \infty, \\
U_{n} ^{\alpha,,\beta,\kappa}f &\rightarrow &f,\text{ \ \ \ a.e., \ \ \ as \ }n\rightarrow\infty ,\text{\ \ \ }\\
V_{n}^{\alpha,\kappa}f &\rightarrow &f,\text{ \ \ \ a.e., \ \ \ as \ }n\rightarrow\infty ,\text{\ \ \ }
\end{eqnarray*}
\end{corollary}
	
\begin{theorem} \label{theorem2}
a) The maximal operator $T^{\ast,\kappa}$ of the summability method (\ref{nor}) with non-decreasing sequence $\{q_{k}:k\geq 0\}$ satisfying the condition 
\begin{equation} \label{Tmeanscond}
\frac{q_{n-1}}{Q_n}= O\left(\frac{1}{n}\right), \ \ \ \text{as} \ \ \ n\to \infty,
\end{equation}
is bounded from the Hardy space $H_{1/2}$ to the space $L_{1/2,\infty}.$
		
Our next result shows that the statement in b) above hold also for non-decreasing sequences and now without any restriction like (\ref{cond1}).
		
b) Let $0<p<1/2.$ For any non-decreasing sequence $\{q_{k}:k\geq 0\},$ there exists a martingale $f\in H_{p},$ such that
\begin{equation*}
\underset{n\in\mathbb{N}}{\sup}\left\Vert T^\kappa_{n}f\right\Vert _{L_{p,\infty}}=\infty.
\end{equation*}		
\end{theorem}
	
A number of special cases of our results are of particular interest and give both well-known and new information. We just give the following examples of such  $T$ means with non-decreasing sequence $\{q_{k}:k\geq 0\}:$
\begin{corollary} \label{cor4}
The maximal operator $B^{\alpha,\beta, *,\kappa}$ is  bounded from the Hardy space $H_{1/2}$ to the space $L_{1/2,\infty}$ but is not bounded from $H_{p}$ to the space $L_{p,\infty},$ when $0<p<1/2.$
\end{corollary}
	
\begin{corollary}\label{cor5}
Let $f\in L_{1}$ and $T^{\kappa}_{n}$ be the $T$ means with \textit{non-decreasing sequence }$\{q_{k}:k\geq 0\}$ and satisfying condition (\ref{Tmeanscond}). Then
\begin{equation*}
T^{\kappa}_{n}f\rightarrow f,\text{ \ \ a.e., \ \ as \ }n\rightarrow \infty .
\end{equation*}
\end{corollary}
	
\begin{corollary}\label{cor6}
Let $f\in L_{1}$. Then $ \ B_{n}^{\alpha,\beta,\kappa}f \rightarrow f,\text{ \ \ \ a.e., \ \ \ as \ }n\rightarrow 	\infty.$
\end{corollary}

\section{Proofs}
\textbf{Proof of Theorem 1.} Let the sequence $\{q_{k}:k\geq 0\}$ be non-increasing. By combining (\ref{2b}) with (\ref{2c}) and using Abel transformation we get that
\begin{eqnarray*}
\left\vert T_{n}^{\kappa}f\right\vert &\leq &\left\vert \frac{1}{Q_{n}}\overset{n-1}{\underset{j=1}{\sum }}q_{j}S_{j}^{\kappa}f\right\vert \\
&\leq &\frac{1}{Q_{n}}\left( \overset{n-2}{\underset{j=1}{\sum }}\left\vert q_{j}-q_{j+1}\right\vert j\left\vert \sigma^{\kappa} _{j}f\right\vert +q_{n-1}(n-1)\left\vert \sigma _{n}^{\kappa}f\right\vert \right) \\
&\leq &\frac{1}{Q_{n}}\left( \overset{n-2}{\underset{j=1}{\sum }}\left( q_{j}-q_{j+1}\right) j+q_{n-1}(n-1)\right) \sigma ^{\ast,\kappa}f\leq \sigma ^{\ast,\kappa}f
\end{eqnarray*}
so that
\begin{equation}
T^{\ast,\kappa }f\leq \sigma ^{\ast,\kappa}f.  \label{12aaaa}
\end{equation}
	
If we apply (\ref{12aaaa}) and Theorem W1 we can conclude that the maximal operators  $T^{\ast,\kappa }$ of all $T$ means with non-increasing sequence $\{q_{k}:k\geq 0\},$ are bounded from the Hardy space $H_{1/2}$ to the space $weak-L_{1/2}.$ 
	
On the other hand by using proposition \ref{prop2} we also have weak $(1,1)$ type inequality and by well-known density argument due to Marcinkiewicz and Zygmund \cite{mz} we have $T_{n}^\kappa f\rightarrow f,$ a.e., for all $f\in L_1.$

Let $0<p<1/2$ and $\left\{ \alpha _{k}:k\in\mathbb{N}\right\} $ be an increasing sequence of positive integers such that:\qquad
\begin{equation}
\sum_{k=0}^{\infty }1/\alpha _{k}^{p}<\infty ,  \label{2}
\end{equation}
	
\begin{equation}
\sum_{\eta =0}^{k-1}\frac{2^{\alpha _{\eta }/p}}{\alpha _{\eta }}<\frac{2^{\alpha _{k}/p-1}}{2\alpha _{k}},  \label{3}
\end{equation}
	
\begin{equation}
\frac{2^{\alpha_{k-1}(1/p-1)}}{\alpha_{k-1}}< \frac{2^{\alpha_k(1/p-1)-4}}{\alpha _{k}}.  \label{4}
\end{equation}
	
We note that such an increasing sequence $\left\{ \alpha _{k}:k\in\mathbb{N}\right\} $ which satisfies conditions (\ref{2}-\ref{4}) can be constructed.
	
Let 
\begin{equation}
f^{\left( A\right) }=\sum_{\left\{ k;\text{ }\lambda _{k}<A\right\} }\lambda_{k}a_{k},  \label{55}
\end{equation}%
where	
\begin{equation*} \label{77}
\lambda _{k}=\frac{1}{\alpha _{k}}  \ \ \ \  \text{and} \ \ \ \ a_k=2^{\alpha_k(1/p-1)}\left(D_{2^{\alpha_k+1}}-D_{2^{\alpha_k}}\right).\end{equation*}
	
By using Proposition \ref{prop1}, it is easy to see that the martingale $\,f=\left( f^{\left( 1\right)},f^{\left( 2\right) }...f^{\left( A\right) }...\right) \in H_{1/2}.$ Moreover, it is easy to show that
\begin{equation} \label{6}
\widehat{f}(j)=\left\{
\begin{array}{ll}
\frac{2^{\alpha _{k}{(1/p-1)}}}{\alpha _{k}},\, & \text{if \thinspace\thinspace }j\in \left\{ 2^{\alpha _{k}},...,2^{\alpha _{k}+1}-1\right\},
\text{ }k=0,1,2..., \\
0, & \text{if \ \thinspace \thinspace \thinspace }j\notin\bigcup\limits_{k=1}^{\infty }\left\{ 2^{\alpha_k},...,2^{\alpha
_{k}+1}-1\right\} 
\end{array}
\right.  
\end{equation}
We can write
\begin{equation}\label{155aba}
T^{\kappa}_{2^{\alpha _{k}}+2}f=\frac{1}{Q_{2^{\alpha _k}+2}}\sum_{j=0}^{2^{\alpha_k}}q_{j}S^{\kappa}_{j}f+\frac{q_{2^{\alpha_k}+1}}{Q_{2^{\alpha _{k}}+2}}S^{\kappa}_{2^{\alpha_{k}}+1}f:=I+II. 
\end{equation}
	
Let $2^{\alpha_s}\leq j\leq 2^{\alpha _{s}+1},$ where $s=0,...,k-1.$
Moreover,
\begin{equation*}
\left\vert D^\kappa_j-D_{2^{\alpha_s}}\right\vert \leq j-2^{\alpha_s}\leq 2^{\alpha_s},\text{ \ \ }\left(s\in\mathbb{N}\right)
\end{equation*}%
so that, according to (\ref{Dir}) and (\ref{6}), we have that
\begin{eqnarray}
&&\left\vert S^\kappa_{j}f\right\vert=\left\vert \sum_{v=0}^{2^{\alpha _{s-1}+1}-1}\widehat{f}(v)\kappa_{v}+\sum_{v=2^{\alpha _{s}}}^{j-1}\widehat{f}(v)\kappa_{v}\right\vert  \label{8} \\
&\leq &\left\vert \sum_{\eta =0}^{s-1}\sum_{v=2^{\alpha_{\eta
}}}^{2^{\alpha_{\eta}+1}-1}\frac{2^{\alpha _{\eta }(1/p-1)}}{\alpha_{\eta }}\kappa_{v}\right\vert +\frac{2^{\alpha _{s}(1/p-1)}}{\alpha_{s}}\left\vert \left( D^\kappa_{j}-D_{2^{\alpha_s}}\right) \right\vert  \notag \\
&=&\left\vert \sum_{\eta =0}^{s-1}\frac{2^{\alpha _{\eta}(1/p-1)}}{\alpha_{\eta}}\left( D_{2^{\alpha _{\eta }+1}}-D_{2^{\alpha_{\eta}}}\right) \right\vert +\frac{2^{\alpha_s(1/p-1)}}{\alpha _{s}}\left\vert \left( D^\kappa_j-D_{2^{\alpha_s}}\right) \right\vert  \notag \\
&\leq & \sum_{\eta =0}^{s-1}\frac{2^{\alpha_{\eta }/p}}{\alpha_{\eta}}+\frac{2^{\alpha _s/p}}{\alpha_s}\leq \frac{2^{\alpha_{s-1}/p+1}} {\alpha_{s-1}}+\frac{2^{\alpha_s/p}}{\alpha_s}\leq \frac{ 2^{\alpha _{k-1}/p+1}}{\alpha_{k-1}}.  \notag
\end{eqnarray}
	
Let $2^{\alpha_{s-1}+1}+1\leq $ $j\leq 2^{\alpha _{s}},$ where $s=1,...,k.$ Analogously to (\ref{8}) we can prove that
\begin{eqnarray*}
&&\left\vert S^{\kappa}_jf\right\vert =\left\vert \sum_{v=0}^{2^{\alpha _{s-1}+1}-1}\widehat{f}(v)\psi _{v}\right\vert =\left\vert \sum_{\eta=0}^{s-1}\sum_{v=2^{\alpha _{\eta }}}^{2^{\alpha _{\eta }+1}-1}\frac{2^{\alpha _{\eta }(1/p-1)}}{\alpha _{\eta }}\kappa_{v}\right\vert \\
&=&\left\vert \sum_{\eta =0}^{s-1}\frac{M_{\alpha _{\eta }}^{1/p-1}}{\alpha_{\eta }}\left( D_{2^{\alpha _{\eta}+1}}-D_{2^{\alpha _{\eta}}}\right) \right\vert \leq \frac{ 2^{\alpha _{s-1}/p+1}}{\alpha_{s-1}}\leq \frac{ 2^{\alpha _{k-1}/p+1}}{\alpha_{k-1}}.
\end{eqnarray*}
Hence,
\begin{equation}
\left\vert I\right\vert \leq \frac{1}{Q_{2^{\alpha _{k}}+2}}\sum_{j=0}^{2^{\alpha_k}}q_j\left\vert S^{\kappa}_{j}f\right\vert \leq \frac{ 2^{\alpha _{k-1}/p+1}}{\alpha _{k-1}}\frac{1}{Q_{M_{\alpha_{k}}+1}}\sum_{j=0}^{2^{\alpha_{k}}}q_{j}\leq \frac{ 2^{\alpha_{k-1}/p+1}}{\alpha _{k-1}}.  \label{10}
\end{equation}
If we now apply (\ref{6}) and (\ref{8}) we get that
\begin{eqnarray}\label{100} 
\left\vert II\right\vert&=&\frac{q_{2^{\alpha_k}+1}} {Q_{2^{\alpha_{k}}+2}}\left\vert \frac{2^{\alpha_{k}(1/p-1)}}{\alpha _{k}}\kappa _{2^{\alpha_k}}+S_{2^{\alpha_k}}f\right\vert  \\
&=&\frac{q_{2^{\alpha_{k}}+1}}{Q_{2^{n_k}+2}} \left\vert \frac{2^{\alpha _{k}(1/p-1)}}{\alpha _{k}}\kappa_{2^{\alpha _{k}}}+S_{2^{\alpha _{k-1}+1}}f\right\vert\notag \\
&\geq&\frac{q_{2^{\alpha_{k}}+1}}{Q_{2^{n_{k}}+2}}\left( \left\vert \frac{2^{\alpha_{k}(1/p-1)}}{\alpha_{k}}\kappa_{2^{\alpha_k}}\right\vert -\left\vert S_{2^{\alpha _{k-1}+1}}f\right\vert \right)  \notag \\
&\geq &\frac{q_{2^{\alpha _{k}}+1}}{Q_{2^{\alpha _k}+2}}\left( \frac{2^{\alpha_{k}(1/p-1)}}{\alpha _k}-\frac{2^{\alpha_{k-1}/p+2}}{\alpha _{k-1}}\right) \notag\\ 
&\geq& \frac{q_{2^{\alpha _{k}}+1}}{Q_{2^{\alpha _{k}}+2}}\frac{2^{\alpha_{k}(1/p-1)-2}}{\alpha _{k}}.  \notag
\end{eqnarray}
	
Without lost the generality we may assume that $c=1$ in (\ref{cond1}). By combining (\ref{10}) and (\ref{100}) we get
\begin{eqnarray}\label{155aba2}
\left\vert T^{\kappa}_{2^{\alpha _{k}}+2}f\right\vert 
&\geq &\left\vert II\right\vert-\left\vert I\right\vert \\ \notag
&\geq &\frac{q_{2^{\alpha _{k}}+1}}{Q_{2^{\alpha_k}+2}}\frac{2^{\alpha _{k}(1/p-1)-2}}{\alpha_k}-\frac{2^{\alpha _{k-1}/p+1}}{\alpha_{k-1}} \\
&\geq &\frac{2^{\alpha _{k}(1/p-2)}}{4\alpha _{k}}-\frac{2^{\alpha_{k-1}/p+1}}{\alpha _{k-1}}\geq \frac{2^{\alpha _{k}(1/p-2)}}{16\alpha _{k}}. \notag
\end{eqnarray}
On the other hand,
\begin{equation}
\mu \left\{ x\in G_{m}:\left\vert T^{\kappa}_{2^{\alpha_k}+2}f\left( x\right)\right\vert \geq \frac{2^{{\alpha _{k}}{(1/p-2)}}}{16\alpha _{k}}\right\} =\mu\left( G_{m}\right) =1.  \label{88}
\end{equation}
Let $0<p<1/2.$ Then
\begin{eqnarray}
&&\frac{2^{{\alpha _{k}}{(1/p-2)}}}{16\alpha _{k}}\cdot \left(\mu \left\{ x\in G_{m}:\left\vert T^{\kappa}_{M_{\alpha _{k}}+2}f\left( x\right) \right\vert \geq \frac{2^{{\alpha _{k}}{(1/p-2)}}}{16\alpha _{k}} \right\} \right)^{1/p} \label{99} \\
&=&\frac{2^{{\alpha _{k}}{(1/p-2)}}}{16\alpha_k}\rightarrow \infty ,\text{ \	as }k\rightarrow \infty .  \notag
\end{eqnarray}
The proof is complete.
	
\textbf{Proof of Corollary \ref{cor0}.} Since $R^\kappa_n,	U_n^{\alpha,\kappa}$ and  $V_{n}^{\alpha,\kappa}$ are  the $T$ means with non-increasing sequence $\{q_{k}:k\geq 0\},$ then the proof of this corollary is direct consequence of Theorem \ref{theorem1}.
	
\textbf{Proof of Corollary \ref{cor2}.} According to Theorem 1 a) and Proposition \ref{prop2} we also have weak $(1,1)$ type inequality and by well-known density argument due to Marcinkiewicz and Zygmund \cite{mz} we have $T^\kappa_{n}f\rightarrow f,$ a.e., for all $f\in L_1.$ Which follows proof of Corollary \ref{cor2}.
	
\textbf{Proof of Corollary \ref{cor1}.} Since $R^\kappa_{n},	U_{n}^{\alpha,\kappa}$ and  $V_{n}^{\alpha,\kappa}$ are  the $T$ means with non-increasing sequence $\{q_{k}:k\geq 0\},$ then the proof of this corollary is direct consequence of Corollary \ref{cor2}.
	
\textbf{Proof of Theorem 2.} Let the sequence $\{q_{k}:k\geq 0\}$ be non-decreasing. By combining (\ref{2b}) with (\ref{2c}) and using Abel transformation we get that

\begin{eqnarray*}
\left\vert T_{n}f\right\vert &\leq &\left\vert \frac{1}{Q_{n}}\overset{n-1}{\underset{j=1}{\sum }}q_{j}S_{j}f\right\vert \\
&\leq &\frac{1}{Q_{n}}\left( \overset{n-2}{\underset{j=1}{\sum }}\left\vert q_{j}-q_{j+1}\right\vert j\left\vert \sigma _{j} f\right\vert+q_{n-1}(n-1)\left\vert \sigma _{n}f\right\vert \right) \\
&\leq &\frac{1}{Q_{n}}\left( \overset{n-2}{\underset{j=1}{\sum }}-\left(q_{j}-q_{j+1}\right) j-q_{n-1}(n-1)+2q_{n-1}(n-1)\right)\sigma ^{\ast }f \\
&\leq &\frac{1}{Q_{n}}\left(2q_{n-1}(n-1)-Q_{n}\right)\sigma ^{\ast }f
\leq c\sigma ^{\ast}f
\end{eqnarray*}
so that
\begin{equation}
T^{\ast }f\leq c \sigma ^{\ast }f.  \label{12aaaaa}
\end{equation}
	
If we apply (\ref{12aaaaa}) and Theorem W1 we can conclude that the maximal operators $T^{\ast }$ is bounded from the Hardy space $H_{1/2}$ to the space $weak-L_{1/2}.$ According to Proposition \ref{prop2} we can conclude that $T^{\ast }$ has weak type-(1,1) and by well-known density argument due to Marcinkiewicz and Zygmund \cite{mz} we also have $T_{n}f\rightarrow f,$ a.e..
	
To prove part b) of theorem 2 we use the martingale, defined by (\ref{55}) where $\alpha_k$ satisfy conditions (\ref{2}-\ref{4}). It is easy to show that for every non-increasing sequence $\{q_{k}:k\geq 0\}$ it automatically holds that
$
q_{2^{\alpha_k+1}}/Q_{2^{\alpha _{k}+2}}\geq c/2^{\alpha_k}.
$
According to (\ref{155aba}-\ref{155aba2}) we can conclude that
\begin{equation*}
\left\vert T^\kappa_{2^{\alpha_k}+2}f\right\vert \geq \left\vert II\right\vert-\left\vert I\right\vert \geq \frac{2^{\alpha _{k}(1/p-2)}}{16\alpha _{k}}.
\end{equation*}

Analogously to\ (\ref{88}) and (\ref{99}) we then get that 
$
\sup_{k\in \mathbb{N}}\left\Vert T^{\kappa}_{2^{\alpha _{k}}+2}f\right\Vert _{L_{p,\infty}}=\infty 
$
and the proof is complete.
	
\textbf{Proof of Corollary \ref{cor4}.} Since $B^{\alpha,\beta, *,\kappa}$ are  the $T$ means with non-decreasing sequence $\{q_{k}:k\geq 0\},$ then the proof of this corollary is direct consequence of Theorem \ref{theorem2}.
	
\textbf{Proof of Corollary \ref{cor5}.} According to Proposition \ref{prop2} we can conclude that $T^{\ast, \kappa}$ has weak type-(1,1) and by well-known density argument due to Marcinkiewicz and Zygmund \cite{mz} we also have $T^\kappa_{n}f\rightarrow f,$ a.e.. Which follows proof of Corollary \ref{cor5}.

\textbf{Proof of Corollary \ref{cor6}.} Since $B^{\alpha,\beta,*,\kappa}$ are  the $T$ means with non-decreasing sequence $\{q_{k}:k\geq 0\},$ then the proof of this corollary is direct consequence of Corollary  \ref{cor5}.

\end{document}